\documentclass[11pt]{article}
\usepackage[utf8]{inputenc}
\usepackage[english]{babel}

\usepackage{amsfonts}
\usepackage{amsmath,amsthm,amscd,amssymb,mathrsfs,setspace}
\usepackage{mathtools}
\usepackage{cite}
\usepackage{latexsym,epsf,epsfig}
\usepackage{color}
\usepackage[a4paper,top=1in,bottom=1in,left=1in,right=1in,marginparwidth=0.5in]{geometry}
\usepackage[colorlinks,citecolor=lccx]{hyperref}
\usepackage{comment}
\usepackage{comment}
\usepackage[nameinlink,capitalise]{cleveref}
\usepackage[mathscr]{euscript}
\usepackage{dutchcal}
\usepackage{algorithm}
\usepackage{algorithmicx}
\usepackage{algpseudocode}
\usepackage[table]{xcolor}
\usepackage{soul}
\usepackage{subcaption}
\usepackage{siunitx}
\usepackage{tikz}
\usetikzlibrary{decorations.pathreplacing,calligraphy}
\usepackage{comment}
\usepackage{booktabs}
\usepackage{multirow}
\usepackage{adjustbox}
\usepackage{authblk}
\usepackage{comment}
\usepackage{enumitem}

\definecolor{lccx}{HTML}{92268F}

\usepackage{amssymb}
\usepackage{xcolor}
\usepackage{verbatim}

\newcommand{\cX}{\mathcal X}
\newcommand{\cY}{\mathcal Y}
\newcommand{\cZ}{\mathcal Z}
\renewcommand{\Re}{\rm Re}
\newcommand{\ran}{\rm ran}
\usepackage{graphics} 
\usepackage{epsfig} 
\usepackage{mathptmx} 
\usepackage{times} 
\usepackage{amsmath}
\usepackage{amssymb} 

\newtheorem{theorem}{Theorem}
\newtheorem{lemma}[theorem]{Lemma}
\newtheorem{proposition}[theorem]{Proposition}
\newtheorem{cor}[theorem]{Corollary}
\newtheorem{defn}[theorem]{Definition}
\newtheorem{example}[theorem]{Example}

\newtheorem{rem}[theorem]{Remark}
\usepackage{bm} 

          \renewcommand{\Re}{\operatorname{Re}}

\title{\LARGE \bf
On Solvability of  Irregular Infinite-Dimensional Differential-Algebraic Equations
}

\author[1]{Ala' Alalabi}
\author[1]{Lorena Bociu}
\author[2]{Birgit Jacob}
\author[3]{Kirsten Morris}

\affil[1]{Department of Mathematics,
North Carolina State University,
Raleigh, NC, USA\\
\texttt{awalalab@ncsu.edu},
\texttt{lvbociu@ncsu.edu}}

\affil[2]{Department of Mathematics,
University of Wuppertal,
42119 Wuppertal, Germany\\
\texttt{bjacob@uni-wuppertal.de}}

\affil[3]{Department of Applied Mathematics,
University of Waterloo,
Waterloo, ON N2L 3G1, Canada\\
\texttt{kmorris@uwaterloo.ca}}

\date{}

\begin{document}

\maketitle
\makeatletter
\begingroup
\renewcommand{\@makefnmark}{}
\footnotemark
\footnotetext{The work of Lorena Bociu was  partially supported by NSF DMS-2510261. Birgit Jacob  and  Kirsten Morris  acknowledge funding by the Deutsche Forschungsgemeinschaft (DFG, German Research Foundation) – Project-ID 531152215 – CRC 1701 and by an NSERC Alliance International Grant from the National Science and Engineering Research Council of Canada. \color{black}}
\endgroup
\makeatother

\begin{abstract}
This paper investigates the solvability of infinite-dimensional
differential-algebraic equations (DAEs) on Hilbert spaces.
We focus on irregular systems for which the operator pencil is
not surjective onto the full state space, but only onto a subspace
determined by the algebraic constraints.
We establish a Lumer-Phillips-type generation result for this
class of systems, showing that dissipativity together with a
partial surjectivity condition yields a contraction semigroup
on the dynamic subspace.
The resulting wellposedness conditions are compared with
existing solvability results.
We further extend the analysis to systems in which hidden
algebraic constraints restrict the dynamics to a proper closed
subspace of the dynamical subspace.
Under an additional structural condition, we show that the
resulting constrained system is $E$-radial on suitably chosen
spaces.
We also demonstrate that semigroup generation does not,
in general, imply $E$-radiality on the original state spaces.
This provides a Weierstra\ss-type decomposition useful
for controller design.
The framework is illustrated by some examples.
\end{abstract}
\newcommand{\keywords}[1]{%
\par\medskip
\noindent\textbf{Keywords:} #1
\par\medskip}

\section{Introduction}

Linear infinite-dimensional differential-algebraic equations (DAEs) arise naturally in systems subject to algebraic or physical constraints. Examples include electrical networks containing semiconductor devices \cite{ali2003elliptic}, and coupled partial differential equation (PDE) systems in which one component is governed by an equilibrium relation. Applications of the latter include electrochemical models of lithium-ion cells \cite{wu2006well,consiglieri2020weak}, biological transport networks \cite{li2020global}, piezoelectric beams with quasi-static magnetic effects \cite{morris2014modeling}, and poroelastic systems \cite{biot1941general,showalter2000diffusion,ARMA,bociu2022weak}.

Infinite-dimensional DAEs can be written in the abstract form
\begin{align}
    \frac{d}{dt}Ex(t)
    &= Ax(t), \qquad t\in(0,T],
    \qquad (Ex)(0)=z_i,
    \label{1}
\end{align}
where $E\in\mathcal{L}(\cX,\cZ)$ and $A:D(A)\subset\cX\to\cZ$ is a closed, densely defined operator. Here, $\cX$ and $\cZ$ are complex separable Hilbert spaces and $z_i\in\cZ$ is the initial condition. If $E$ is boundedly invertible, then \eqref{1} reduces to a classical evolution equation. We are therefore most interested in the degenerate case in which $E$ has a nontrivial kernel. Equations of this type are also referred to as implicit evolution equations \cite{showalter2010hilbert,showalter1996monotone}, since the evolution cannot, in general, be written explicitly as an equation for $\dot{x}(t)$. Such systems are also referred to as partial differential-algebraic equations (PDAEs).

The wellposedness of infinite-dimensional DAEs has been studied from several operator-theoretic perspectives; see, for example,
\cite{showalter2010hilbert,jacob2022solvability,sviridyuk2003linear, zwart2024abstract,thaller1999semigroup,favini1998degenerate}.
Thaller and Thaller \cite{thaller1999semigroup} studied solvability in Hilbert spaces using the orthogonal decompositions
$$     \cX=\ker E\oplus\overline{\operatorname{ran}E^*},  \qquad
    \cZ=\ker E^*\oplus\overline{\operatorname{ran}E}. $$
Sviridyuk and Fedorov \cite{sviridyuk2003linear} developed a different approach based on $(E,p)$-radiality and Hille-Yosida-type resolvent estimates. Related radiality conditions were subsequently relaxed in \cite{jacob2022solvability,erbay2024weierstrass}. More recently, Zwart and Mehrmann \cite{zwart2024abstract} studied solvability and state-space reduction for a class of infinite-dimensional DAEs possessing a dissipative Hamiltonian structure.

Many existing solvability results are formulated under regularity
assumptions on the operator pencil $(E,A)$. In particular, they require suitable injectivity and surjectivity properties of
$    sE-A:D(A)\to\cZ $
for $s$ in an appropriate resolvent set. In some infinite-dimensional systems, however, surjectivity onto the full space $\cZ$ is too restrictive. In these systems, it is possible that the dynamics can be defined on  a smaller subspace.  We refer to such systems as \emph{irregular}. A class of irregular implicit evolution equations was considered by Showalter
\cite[Chapter 5]{showalter2010hilbert}, where dissipativity is combined with a surjectivity condition on $\operatorname{ran}E$ rather than on the full space $\cZ$.

The main objective of this paper is to develop a wellposedness framework for irregular dissipative infinite-dimensional DAEs. We show that dissipativity, together with an appropriate partial surjectivity condition, gives rise to a reduced operator that generates a contraction semigroup on the dynamic subspace. This yields a Lumer-Phillips-type generation result without requiring surjectivity of the operator pencil on the entire state-space. We compare the resulting conditions with existing solvability results for infinite-dimensional DAEs given in \cite{showalter2010hilbert} and \cite{zwart2024abstract}.
The results here generalize those in \cite{showalter2010hilbert,zwart2024abstract} and provide a unified presentation. 

We also investigate the relation between this reduced semigroup framework and $E$-radiality. We show that generation of a contraction semigroup on the dynamic subspace does not, in general, imply $E$-radiality on the original spaces, since the required resolvent estimates may fail. We identify an additional structural condition under which $E$-radiality can be recovered on appropriately chosen spaces. The resulting state-space decomposition is analogous to the Weierstra\ss{} form for finite-dimensional DAEs. Such decompositions are particularly relevant for the analysis and design of controllers  for infinite-dimensional DAEs; see, for example,  \cite{alalabi2024finite}.

We also show how the abstract conditions arise naturally in a variational setting based on sesquilinear forms. In particular, coercivity yields the solvability and resolvent estimates required for $E$-radiality. Finally, the theoretical framework is illustrated on a class of coupled systems arising in physically relevant models \cite{morris2025boundary,bociu2025analysis,chaturvedi2010algorithms}.
Our assumptions weaken those imposed in \cite[Section IV]{jacob2022solvability}.

The paper is organized as follows. Section \ref{section2} develops the wellposedness theory for irregular dissipative infinite-dimensional DAEs, compares the resulting conditions with existing solvability results, and investigates their relation to $E$-radiality. Section \ref{section1} recalls the $(E,p)$-radiality framework and the associated state-space decomposition. The subsequent sections develop the variational formulation and extend the generation and $E$-radiality results to systems whose
dynamics evolve on a constrained subspace. Finally, the abstract theory is illustrated on a class of coupled systems.

A preliminary version of some of the results in this paper was presented at the {\em European Control Conference} \cite{AlalabiECC}.
In this paper the presentation is changed, along with several proofs. A more complete comparison with \cite{showalter2010hilbert,zwart2024abstract} is provided, and  the constrained subspace approach in Section \ref{section4} is new. \\

\noindent\textit{Notation.}
For Hilbert spaces $X$ and $Y$, $\mathcal{L}(X,Y)$ denotes the space of bounded linear operators from $X$ to $Y$. The kernel and range of an operator $A$ are denoted by $\ker A$ and $\operatorname{ran}A$, respectively. We denote the open right half-plane by
$ \mathbb{C}_+  :=   \{s\in\mathbb{C}:\operatorname{Re}s>0\}.$

\section{wellposedness of Irregular Systems} \label{section2}
For Hilbert spaces $\cX, \cZ,$   let $E\in\mathcal{L}(\cX,\cZ)$ and $A:D(A)\subset \cX \to \cZ$ be a closed, densely defined operator.  We define the  \emph{$E$-resolvent} with respect to $A$  as $(sE - A)^{-1},$ and the resolvent set  $$\rho(E,A) = \left\{ s \in \mathbb{C} : (sE -A)^{-1} \in \mathcal{L}(\mathcal{Z}, \mathcal{X}) \right\}. $$

 As in  the finite-dimensional settings \cite{kunkel1997linear}, we use the following terminology.  
\begin{defn}
The pair $(E,A)$ is \emph{regular} if there exists $s \in \mathbb{C}$ such that the operator $s E - A: D(A) \to \mathcal{Z}$ is boundedly invertible. If no such $s$ exists, the pair $(E,A)$ is \emph{irregular}.
\end{defn}

For irregular systems we seek wellposedness on an appropriate subspace of the original state space.
 The general system \eqref{1} is specialized to systems possessing an additional dissipative Hamiltonian structure. More precisely, we consider
\begin{align}
    \frac{d}{dt} (Ex(t)) =& AQ x(t), \quad t\ge 0, \; \; \; (Ex)(0) = z_i, \label{system1}
\end{align}
where $E\in\mathcal{L}(\cX,\cZ),$ $A:D(A)\subset \cZ\to\cZ$ is densely defined, and $Q\in\mathcal{L}(\cX,\cZ)$ is an isomorphism. If, in addition, $\operatorname{ran} E$ is closed in $\cZ$,
and $E^\ast Q = Q^\ast E \geq 0$, then
$Q^\ast E=E^\ast Q$ also has closed range.
Moreover, $Q^\ast E$ is self-adjoint, and, since $Q^\ast$
is invertible, $\ker(Q^\ast E)=\ker(E)$. 
Thus, we have the orthogonal decomposition
$    \cX=\operatorname{ran}(Q^\ast E)\oplus\ker(E).$
We denote
\begin{align} \cX_1&:=\operatorname{ran}(Q^\ast E),  \qquad \cX_2:=\cX_1^\perp=\ker(E).  \label{orthogonal-decomposition}
\end{align}
 We define $\widetilde E\in L(\cX_1)$ by $\widetilde Ex:= Q^\ast Ex, \; \;  x\in \cX_1.$  The restriction $\widetilde E$ is injective.  Furthermore, $\operatorname{ran}\widetilde E =\operatorname{ran}(Q^*E) =\mathcal X_1, $ because $Q^*E$ vanishes on $\mathcal X_2$. Thus $\widetilde E$ is bijective on $\mathcal X_1$. By the bounded inverse theorem, $\widetilde E^{-1}\in\mathcal L(\mathcal X_1).$   Further, $ EQ^{-1} = Q^{-\ast}Q^{\ast} EQ^{-1}  $ is a self-adjoint operator on $\cZ$ with $ EQ^{-1}=  Q^{-\ast}E^\ast \ge 0$.
Thus, we can write  $\cZ$ as the direct orthogonal sum
$ \cZ=  \mathrm{ran} \, E \oplus \mathrm{ker}(EQ^{-1}),$
and we denote $$ \cZ_1:=\mathrm{ran} \, E , \qquad   \cZ_2:=\cZ_1^\bot=\mathrm{ker}(E Q^{-1}). $$ We note that $Q\cX_2 = \cZ_2 $.

We now recall the following well-established result.
\begin{proposition}(Lumer-Phillips Theorem, \cite[eg., Theorem 6.1.7]{JZ})\label{prop2b}
 Let $B$ be a linear operator with domain
$D(B)$ on a Hilbert space $\cY$. Then $B$ is the generator of a contraction
semigroup  on $\cY$ if and only if $B$ is dissipative and  $I-B$ is surjective.
\end{proposition}

It will now be shown, using the Lumer-Phillips Theorem, that dissipativity on a subspace $\cX_1$, together with a range condition restricted to $\operatorname{ran}E$, yields a contraction semigroup. Regularity, that is, surjectivity of the operator pencil onto the entire space $\cZ,$ is not assumed.

\begin{theorem}\label{thm1b} Consider system \eqref{system1}. Suppose that
$\operatorname{ran} E$ is closed in $\mathcal{Z}$ and
$E^\ast Q=Q^\ast E\geq 0$.
Recall the spaces defined in \eqref{orthogonal-decomposition}.
Equip $\cX_1$ with the inner product
$\langle x,y\rangle_{\cX_1}:=
\langle x,\widetilde E y\rangle_{\cX}$
for $x,y\in\cX_1$, and define
$D:=\{x\in D(AQ):Q^*AQx\in\cX_1\}$.  Assume that
\begin{enumerate}
\item $\left.Q^*AQ\right|_D$ is dissipative;
\item $\cX_2\cap D=\{0\};$
\item $\operatorname{ran}E\subseteq\operatorname{ran}(E-AQ)$.
\end{enumerate}
Then the operator
$A_r:D(A_r)\subset\cX_1\to\cX_1$, defined by $A_rx_1:=\widetilde E^{-1}Q^*AQ(x_1+x_2)$,  
with domain $$D(A_r):=\{x_1\in\cX_1:\text{ there exists }x_2\in\cX_2 \text{ such that }x_1+x_2\in D(AQ)\text{ and }Q^*AQ(x_1+x_2)\in\cX_1\},$$
is well-defined and generates a contraction semigroup on $(\cX_1,\langle\cdot,\cdot\rangle_{\cX_1})$.
\end{theorem}
\begin{proof}
We first show that $A_r$ is well-defined. Let $x_1\in\cX_1$ and suppose that $x_2,\widetilde x_2\in\cX_2$ both satisfy the conditions in the definition of $D(A_r)$. Then $y:=x_2-\widetilde x_2\in\cX_2\cap D$ since $y\in D(AQ)$ and $Q^*AQy\in\cX_1$. By the second assumption, $y=0$, and hence $x_2=\widetilde x_2$. Thus, $A_r$ is well-defined.

Next, let $x_1\in D(A_r)$ and choose $x_2\in\cX_2$ as in the definition of $D(A_r)$. Since $Q^*AQ(x_1+x_2)\in\cX_1$ and $x_2\in\cX_2=\cX_1^\perp$, we obtain
\begin{align*}
\Re\langle A_rx_1,x_1\rangle_{\cX_1}
&=\Re\langle Q^*AQ(x_1+x_2),x_1\rangle_{\cX}=\Re\langle Q^*AQ(x_1+x_2),x_1+x_2\rangle_{\cX}
\leq 0,
\end{align*}
where the last inequality follows from the first assumption. Hence, $A_r$ is dissipative.

By the Lumer-Phillips theorem (Proposition \ref{prop2b}), it remains to show that $I-A_r$ is surjective. Let $y_1\in\cX_1$. By the third assumption, there exist $x_1\in\cX_1$ and $x_2\in\cX_2$ such that $x_1+x_2\in D(AQ)$ and $Ey_1=(E-AQ)(x_1+x_2)$. Since $Ex_2=0$,
\[
Q^*AQ(x_1+x_2)=Q^*Ex_1-Q^*Ey_1=\widetilde E(x_1-y_1)\in\cX_1.
\]
Thus, $x_1\in D(A_r)$ and $A_rx_1=x_1-y_1$, so $(I-A_r)x_1=y_1$. Therefore, $I-A_r$ is surjective, and the result follows from the Lumer-Phillips theorem.
\end{proof}

\begin{lemma} \label{showalter-equivalence}
Under the structural assumptions of Theorem~\ref{thm1b},
assume also that $Q=I$, so that
$\cX_1=\operatorname{ran}E$ and $\cX_2=\ker E$. Define  $D: =\{x\in D(A):Ax\in\cX_1\}.$ Then 
\begin{equation} \label{cond_Thm5} \cX_2\cap D = \{ 0 \}  \qquad [\text{Assumption (2) in Theorem \ref{thm1b}}]\end{equation}
if and only if 
    \begin{equation} \label{cond_Show} \cX_2 \cap D \subseteq \ker A \qquad \text{and} \qquad  \cX_2 \cap \ker A = \{0\}.  \qquad [\text{Assumptions in  \cite[Theorem 5.4.1]{showalter2010hilbert} }] \end{equation}
\end{lemma}
\begin{proof}
We first show that \eqref{cond_Show} implies \eqref{cond_Thm5}. Let $x\in\mathcal X_2\cap D$. By the first condition of \eqref{cond_Show}
we have $Ax=0$. Thus, the second condition of \eqref{cond_Show}
implies that $x=0$. Therefore, $\mathcal X_2\cap D=\{0\}$.

Now suppose that \eqref{cond_Thm5} holds. Then the first condition of \eqref{cond_Show} follows immediately.
To verify the second condition of \eqref{cond_Show}, let $x\in\mathcal X_2\cap D(A)$ satisfy
$Ax=0$. Then $Ax\in\mathcal X_1$, and hence $x\in D$. Therefore,
$x\in\mathcal X_2\cap D$, and \eqref{cond_Thm5} yields $x=0$, which completes the proof.
\end{proof}
 
\begin{rem}
Lemma  \ref{showalter-equivalence} shows that, for  the special case \(Q=I\), the kernel condition in
Theorem~\ref{thm1b} is equivalent to the two kernel
conditions in \cite[Theorem 5.4.1]{showalter2010hilbert}. Moreover, the range condition
 in Theorem \ref{thm1b} is equivalent to the surjectivity of $E-A:\mathcal D\to\operatorname{ran}E,
$ where $\mathcal D
=\{x\in D(A):Ax\in\operatorname{ran}E\}.$  However, \cite[Theorem 5.4.1]{showalter2010hilbert} does not require
$\operatorname{ran}E$ to be closed. Instead, the
dynamic space is constructed through the completion
of $\mathcal X/\ker E$ with respect to the seminorm
induced by $E$. Under our additional closed-range
assumption, this space can be identified with
$\mathcal X_1=\operatorname{ran}E$, equipped with the
inner product induced by $E$.

Therefore, Theorem \ref{thm1b} recovers the wellposedness result of Theorem 5.4.1 in \cite{showalter2010hilbert} within this closed-range setting. In addition, unlike \cite[Theorem 5.4.1]{showalter2010hilbert}, Theorem \ref{thm1b} explicitly identifies the dynamic subspace $\cX_1$ and the reduced operator $A_r$, and shows that the semigroup generated by $A_r$ on $\cX_1$ is a contraction semigroup.  Therefore, Theorem \ref{thm1b} provides a direct
representation of the reduced dynamics. It also will
facilitate the  analysis of
$E$-radiality in the next section.
\end{rem}

\begin{rem}
The general setting of Theorem \ref{thm1b} includes, as a special case, the block structure considered in \cite[Theorem 25]{zwart2024abstract},
\begin{align}
E&=
\begin{bmatrix}
E_1&0\\
0&0
\end{bmatrix},
\qquad
Q=
\begin{bmatrix}
Q_1&0\\
0&Q_2
\end{bmatrix},
\qquad
A(x_1,x_2)=
\begin{bmatrix}
A_1(x_1,x_2)\\
A_2(x_1,x_2)
\end{bmatrix},
\label{block-operators}
\end{align}
where $\cX=\cZ=\mathsf X_1\oplus\mathsf X_2$ , also $A:D(A)\subset\cZ\to\cZ$ is densely defined, the operators $E_1,Q_1,$ and $Q_2$ are boundedly invertible and $E_1^\ast Q_1$ is coercive. In this case, the spaces defined in \eqref{orthogonal-decomposition} satisfy $\cX_1=\mathsf X_1\oplus\{0\}$ and $\cX_2=\{0\}\oplus\mathsf X_2$. The wellposedness result in \cite[Theorem 25]{zwart2024abstract} is formulated for regular systems and assumes the dissipativity of $A$ on the full domain $D(A)$. In contrast, Theorem \ref{thm1b} does not require regularity of the operator pencil or the dissipativity on the full domain. Instead, it requires dissipativity of $Q^\ast AQ$ only on the constrained domain $D=\{x\in D(AQ):Q^\ast AQx\in\cX_1\}$, together with the partial range condition $\operatorname{ran}E\subseteq\operatorname{ran}(E-AQ)$. Thus, in the block setting \eqref{block-operators}, the assumptions used in \cite[Theorem 25]{zwart2024abstract} imply the corresponding assumptions of Theorem~\ref{thm1b}. This connection is made precise in the following corollary.
\end{rem}

\begin{cor}\label{comparison-zwart}
Consider the block setting \eqref{block-operators}, and define $D:=\{x\in D(AQ):Q^\ast AQx\in\cX_1\}$. Assume that the hypotheses of \cite[Theorem 25]{zwart2024abstract} hold, namely:
\begin{enumerate}
    \item $A$ is dissipative on $D(A)$;
    \item $\cX_2\cap D=\{0\}$;
    \item the pair $(E,AQ)$ is regular.
\end{enumerate}
Then the assumptions of Theorem \ref{thm1b} are satisfied. Consequently, the reduced operator $A_r$ generates a contraction semigroup on $\cX_1$.
\end{cor}

\begin{proof}
Since $E_1^\ast Q_1$ is self-adjoint and coercive,
$Q^\ast E=E^\ast Q\geq0$.
Moreover, $Q$ is boundedly invertible and
$\operatorname{ran}E=\mathsf X_1\oplus\{0\}$ is closed.
Thus, the structural assumptions of
Theorem \ref{thm1b} are satisfied, with
$\operatorname{ran}(Q^\ast E)=\mathsf X_1\oplus\{0\}$
and $\ker E=\{0\}\oplus\mathsf X_2$.

For $x\in D$, we have $\operatorname{Re}\langle Q^\ast AQx,x\rangle_{\cX}=\operatorname{Re}\langle AQx,Qx\rangle_{\cX}\leq0$, so the first assumption of Theorem \ref{thm1b} holds. The second assumption is precisely $\cX_2\cap D=\{0\}$.

By \cite[Theorem 13]{zwart2024abstract}, the regularity of $(E,AQ)$, together with the dissipative Hamiltonian structure assumed here, implies that $sE-AQ$ is boundedly invertible for every $s\in\mathbb C_+$. In particular, $E-AQ$ is surjective, and hence $\operatorname{ran}E\subseteq\operatorname{ran}(E-AQ)$. Thus, the third assumption of Theorem \ref{thm1b} holds.
\end{proof}

The weaker range condition in Theorem \ref{thm1b} does not, in general, imply regularity of the operator pencil. In particular, the operator pencil need not be surjective onto the full space, as the following example demonstrates.

\begin{example}
Let $\mathcal X_1=\mathcal X_2=\ell^2$ and 
$E=\begin{bmatrix}I&0\\0&0\end{bmatrix}$.
Define $A:D(A)\subseteq \mathcal X_1\times\mathcal X_2\to\ell^2\times\ell^2$ by 
$$D(A):=\ell^2\times\ell^2,\quad A\begin{bmatrix}(x_n)_n\\(y_n)_n\end{bmatrix}
:= \begin{bmatrix} (-x_n)_n\\ (-\frac1n y_n)_n \end{bmatrix}. $$ Taking $Q=I$, we have $E^\ast Q=Q^\ast E=E\geq 0$,
and $\operatorname{ran}E=\ell^2\times\{0\}$ is closed. Clearly, $A$ is bounded and dissipative, and $A\begin{bmatrix}0\\ (y_n)_n\end{bmatrix}\in\ell^2\times\{0\}$
implies $(y_n)_n=0$.  Further, $A_r=-I\in\mathcal L(\ell^2),$ so $A_r$ generates a contraction semigroup on $\ell^2$. Moreover, $\operatorname{ran}E=\ell^2\times\{0\}\subseteq
\operatorname{ran}(E-A)$, so the partial range condition in Theorem \ref{thm1b} is satisfied. However, for $s\in\mathbb C_+$,
$$(sE-A) \begin{bmatrix} (x_n)_n\\ (y_n)_n \end{bmatrix}
= \begin{bmatrix} ((s+1)x_n)_n\\ (\frac1n y_n)_n \end{bmatrix}.$$
Thus, $(sE-A)^{-1}$ is not a bounded operator for any $s$ and hence $(E,A)$ is not regular.
\end{example}

\subsection{Sesquilinear Forms}\label{section3}
The following result shows that, when $Q=I,$ $E$ and $-A$ are  associated with a nonnegative sesquilinear form $m$ and a coercive sesquilinear form $\ell$, respectively, the dissipativity, range, and injectivity conditions
required in the irregular DAE framework of Theorem \ref{thm1b} follow naturally from the variational structure.

\begin{proposition}\label{prop:Showalter-structure}
Let $\cX_m$ be the seminormed space induced by a symmetric,
nonnegative sesquilinear form $m(\cdot,\cdot)$, and let
$E\in\mathcal L(\cX_m,\cX_m^\prime)$ denote the associated bounded operator, defined by $\langle Ex,y\rangle_{\cX_m^\prime,\cX_m}=m(x,y)$.
Suppose that:
\begin{enumerate}
    \item $\cX$ is a Hilbert space, densely and continuously embedded in $\cX_m$, i.e.\ $\cX\hookrightarrow\cX_m$;

    \item $\ell(\cdot,\cdot)$ is a continuous and coercive
    sesquilinear form on $\cX$, and let $-A:\cX\to\cX^\prime$ denote the
    associated operator, defined by
    $\langle -Ax,y\rangle_{\cX^\prime,\cX}=\ell(x,y)$.
\end{enumerate}
Define $D:=\{x\in\cX:Ax\in \cX_m^\prime \}$. Then:
\begin{enumerate}
    \item $A:D\to\cX_m^\prime$ is dissipative;
    \item $E-A:D\to\cX_m^\prime$ is surjective;
    \item $\ker E\cap\ker A=\{0\}$ and $\ker E\cap D=\{0\}$.
\end{enumerate}
\end{proposition}

\begin{proof}
\begin{enumerate}
\item Since $\ell$ is coercive, there exists $c>0$ such that
$\Re\,\ell(x,x)\geq c\|x\|_{\cX}^2$ for all $x\in\cX$.
For $x\in D$, we have $Ax\in\cX_m^\prime$, and hence
$\Re\langle -Ax,x\rangle_{\cX_m^\prime,\cX_m}=\Re\langle -Ax,x\rangle_{\cX^\prime,\cX}=\Re\,\ell(x,x)\geq c\|x\|_{\cX}^2. $
Therefore, $\Re\langle Ax,x\rangle_{\cX_m^\prime,\cX_m}\leq0$, and $A$ is dissipative on $D$.

\item For $\lambda>0$, define $b_\lambda(x,y):=\lambda m(x,y)+\ell(x,y)$ for $x,y\in\cX$. Since $m$ is nonnegative and $\ell$ is coercive, $b_\lambda$ is continuous and coercive on $\cX$. By the Lax-Milgram theorem, for every $w\in\cX_m^\prime$ there exists a unique $x\in\cX$ such that $(\lambda E-A)x=w$. Since $Ex,w\in\cX_m^\prime$, we have $Ax=\lambda Ex-w\in\cX_m^\prime$, and hence $x\in D$. Thus $\lambda E-A:D\to\cX_m^\prime$ is surjective for every $\lambda>0$. In particular, taking $\lambda=1$ gives the surjectivity of $E-A$.

\item By the Lax-Milgram theorem applied to $\ell$, the operator
$-A:\cX\to\cX^\prime$ is an isomorphism. Hence $\ker A=\{0\}$, and
consequently $\ker E\cap\ker A=\{0\}$. Moreover, let $x\in\ker E\cap D$.
    Since $Ex=0$, we have $m(x,x)=0$, and hence
    $\|x\|_{\mathcal{X}_m}=0$.
    Since $Ax\in\mathcal{X}_m'$, it follows that
    $
    \langle Ax,x\rangle_{\mathcal{X}_m',\mathcal{X}_m}=0.
    $
    Therefore,
    \[
    0=\operatorname{Re}\ell(x,x)
    \ge c\|x\|_{\mathcal{X}}^2,
    \]
    which implies $x=0$. Thus
    $\ker E\cap D=\{0\}$.
\end{enumerate}
\end{proof}

The following corollary shows that, under the
assumptions of Proposition~\ref{prop:Showalter-structure} and an additional
closed-range condition, Theorem~\ref{thm1b} applies
with $Q=I$, yielding a contraction semigroup on the
reduced dynamic subspace.

\begin{cor} \label{cor-sesquilinear}Assume the hypotheses of
Proposition \ref{prop:Showalter-structure}. Let $E_{\cX},A_{\cX}\in\mathcal L(\cX)$
be the bounded operators associated with
$m$ and $-\ell$, respectively, through the
Riesz representation theorem on $\cX$, i.e., $$ \langle E_{\cX}x,y\rangle_{\cX}=m(x,y),
 \qquad \langle A_{\cX}x,y\rangle_{\cX}=-\ell(x,y).
$$ Suppose, in addition, that
$\operatorname{ran}E_{\cX}$ is closed. Then the assumptions of
Theorem \ref{thm1b} are satisfied for
$(E_{\cX},A_{\cX})$ with $Q=I$. Then the reduced operator
$A_r$ defined in Theorem \ref{thm1b}
generates a contraction semigroup on $\cX_1:=\operatorname{ran}E_{\cX},
$ equipped with the inner product $ \langle x,y\rangle_{\cX_1}
 : = \langle x,E_{\cX}y\rangle_{\cX}. $
\end{cor}
\begin{proof} Since $m$ is symmetric and nonnegative,
$E_{\cX}$ is self-adjoint and nonnegative.
Together with the closed-range assumption,
this verifies the structural assumptions of
Theorem~\ref{thm1b} with $Q=I$. Moreover, since $\operatorname{ran}E_{\cX}$
is closed, we have $\cX=\operatorname{ran}E_{\cX} \oplus\ker E_{\cX}. $ The constrained domain
$\mathcal D_{\cX} :=\{x\in\cX: A_{\cX}x\in\operatorname{ran}E_{\cX}\}$
coincides with the domain $D$ in
Proposition~\ref{prop:Showalter-structure}.
Thus, the dissipativity, kernel, and range
conditions established in that proposition
yield the corresponding assumptions of
Theorem~\ref{thm1b} for
$(E_{\cX},A_{\cX})$. Since $A_{\cX}\in\mathcal L(\cX)$ is
closed and densely defined, all the
hypotheses of Theorem~\ref{thm1b}
are satisfied, and the conclusion follows. \end{proof}

\begin{rem}
    The assumptions of Proposition \ref{prop:Showalter-structure} coincide with
those of \cite[Thm. 5.4.2]{showalter2010hilbert}.
Proposition~\ref{prop:Showalter-structure} shows that
these variational assumptions imply the dissipativity,
range, and kernel conditions of
\cite[Thm. 5.4.1]{showalter2010hilbert}.
Thus, the operator-theoretic conditions of
\cite[Thm. 5.4.1]{showalter2010hilbert}
are more general than the variational assumptions
of \cite[Thm. 5.4.2]{showalter2010hilbert}.
The latter theorem establishes existence and
uniqueness of solutions to system~\eqref{system1},
assuming that $E$ is induced by a symmetric,
nonnegative sesquilinear form and $-A$ is induced
by a continuous, coercive sesquilinear form on
a densely and continuously embedded Hilbert space.

Furthermore, Corollary~\ref{cor-sesquilinear}
shows that, if the operator $E_{\cX}$ associated
with $m$ on $\cX$ has closed range, the assumptions
of Theorem~\ref{thm1b} are satisfied for
$(E_{\cX},A_{\cX})$ with $Q=I$.
Consequently, the reduced operator $A_r$ generates
a strongly continuous contraction semigroup on
$\cX_1=\operatorname{ran}E_{\cX}$, equipped with
the inner product induced by $E_{\cX}$.
\end{rem} \color{black}

\section{\texorpdfstring{$E$-}{E}radiality theory} 
\label{section1}

For $s\in \rho(E,A),$ we define the right- and the left- $E$ resolvents with respect $A$  as
\begin{align*}
    R^E(s,A)= (sE-A)^{-1} E, \; \; \;  L^E(s,A)= E (sE-A)^{-1} .
\end{align*}
\begin{defn} \label{defn-radiality}
    The operator $A$ is \emph{$E$-radial}, or \emph{$(E,0)$-radial}, if 
    \begin{itemize}
        \item $s\in \rho(E,A)$ for all $s>0,$
        \item $\text{there exists } K>0$ such that for all $n \in \mathbb{N}$ and for all real $s>0$
        \begin{align}
&            \| \left( R^E(s,A)\right)^n\|_{\mathcal{L}(\mathcal{X}
            )} \leq \frac{K}{s^n} , \label{estimate1}\\
            &            \| \left( L^E(s,A)\right)^n\|_{\mathcal{L}(\mathcal{Z}
            )} \leq \frac{K}{s^n} . \label{estimate2} 
        \end{align}
    \end{itemize}
\end{defn}

When $E=I$ in system \eqref{1}, statements \eqref{estimate1} and \eqref{estimate2} coincide, and the Hille-Yosida theorem then yields generation of a $C_{0}$-semigroup.  
The  $(E,p)$-radial framework from \cite{sviridyuk2003linear} reduces to $E$-radiality when $p=0$.

 In parallel to \cite{sviridyuk2003linear}, for $s \in \rho(E,A), $ we define the spaces  
  \begin{subequations}\label{spaces}
      \begin{align}  
    \mathcal{X}^0&= \ker\, (sE - A)^{-1} E, \qquad   \mathcal{Z}^0= \ker\,  E(sE - A)^{-1},  \\
      \mathcal{X}^1&= \overline{\operatorname{ran } (sE - A)^{-1} E}, \qquad   \mathcal{Z}^1= \overline{\operatorname{ran }  E(sE - A)^{-1}}.
\end{align}
  \end{subequations}
These spaces do not depend on the choice of $s$ \cite{sviridyuk2003linear}.
 Note that $ \mathcal{X}^{0} = \ker\,(E),$ and $ \mathcal{X}^0 \cap  \mathcal{X}^1 = \{ 0\}$. Also,   $z \in \mathcal{Z}^{0}$  iff $x=(s E - A)^{-1}z \in \ker\,(E)$.

 Then, since $x \in \ker\,(E)$, it follows that 
 $z=(sE-A)(s E - A)^{-1}z = (sE-A)x = -Ax,$ for some $x \in D(A) \cap \ker\,(E),$ and 
$\mathcal{Z}^{0} = \{ Ax : x \in D(A) \cap \ker\,(E) \}$ \cite{jacob2022solvability}. 
 The following result was obtained in \cite{sviridyuk2003linear}.
\begin{theorem}\cite[Lemma 2.2.1, Lemma 2.2.6, and Corollary 2.5.2]{sviridyuk2003linear}  \label{radiality}  
 If  $A$ is $E-$radial, then
\begin{enumerate}
    \item The spaces $ \mathcal{X}$ and $\mathcal{Z}$ decompose as
$\mathcal{X}= \mathcal{X}^0 \oplus \mathcal{X}^1, \; \; \mathcal{Z}= \mathcal{Z}^0 \oplus \mathcal{Z}^1, $
 with $ \mathcal{X}^0$, $ \mathcal{X}^1$, $\mathcal{Z}^0$ and $\mathcal{Z}^1$   given in \eqref{spaces}.
\item  There exist  bounded linear operators    
 \begin{subequations}
     \begin{align}
  & \mathbf{P}: \mathcal{X} \rightarrow \mathcal{X}^1, \ \mathbf{P} x : =  \lim_{ s \to \infty} s(sE - A)^{-1}E x, \\
&       \mathbf{Q}: \mathcal{Z}  \rightarrow \mathcal{Z}^1, \    \mathbf{Q} z : =  \lim_{ s \to \infty} sE(sE - A)^{-1} z ,  \label{projection2}
\end{align}
 \end{subequations}
 where $  \mathbf{P} $ is a projection onto $ \mathcal{X}^1$ with $\ker\, \mathbf{P} = \mathcal{X}_0, \; \operatorname{ran }  \mathbf{P} = \mathcal{X}^1$, and  $ \mathbf{Q}$ is a projection onto $\mathcal{Z}^1$ with $ \ker\,  \mathbf{Q} = \mathcal{Z}^0, \; \operatorname{ran } \mathbf{Q}= \mathcal{Z}^1$. 

\item Define the restrictions of $E$ and $A$ as 
 \begin{subequations}
\begin{align}
E_0 &:= E|_{ \mathcal{X}^0}, \; A_0 := A|_{D(A_0)}, \; D(A_0) =  D(A) \cap \mathcal{X}^0, \\
E_1 &:= E|_{ \mathcal{X}^1}, \; 
A_1 := A|_{D(A_1)}, \; D(A_1) =  D(A) \cap \mathcal{X}^1 .
\end{align}
 \end{subequations}
For any $x = x_1 + x_0 \in \mathcal{X}$ and for $t>0$, system \eqref{1} can be written in Weierstra\text{\ss} form as
\begin{subequations} \label{decomposition}
    \begin{align}
   \frac{d}{dt}  E_1 x_1(t)  &= A_1 x_1(t) , \label{dynamical subsystem}\\
    0  &= A_0 x_0(t) , \label{algebraic subsystem}\end{align}
    \end{subequations}
where $E_1 \in \mathcal{L}( \mathcal{X}^1, \mathcal{Z}^1),$ $A_1 : D(A_1) \subset \mathcal{X}^1 \to \mathcal{Z}^1$ is densely defined and closed, $A_0 : D(A_0) \subset \mathcal{X}^0 \to \mathcal{Z}^0$ is densely defined, closed, and boundedly invertible.
\end{enumerate}
\end{theorem}

The projections $\mathbf{P}, \mathbf{Q}$ are in general non-orthogonal.

\begin{theorem} \cite[Theorem  1]{jacob2022solvability}
 If $A-\alpha E$ is $E-$radial and $\operatorname{ran} \,E$ is closed, then the reduced operator $E_1^{-1}A_1$, with domain {$D(E_1^{-1}A_1) = D(A_1),$} generates a $C_0$-semigroup $(T(t))_{t\ge 0}$ on $\mathcal{X}^1$ with
bound $\|T(t)\|\le Ke^{\alpha t}$.  In addition, the component of the solution on $\mathcal{X}^0$ remains identically zero.
\end{theorem}

The wellposedness conditions of Theorem \ref{thm1b}, although sufficient for contraction semigroup generation on the dynamic subspace, do not in general imply $E$-radiality. The following example illustrates this distinction.

\begin{example} \label{example-no-rad}
Let $Z=L^2(0,1)$ and define $A_d=\frac{d^2}{d\xi^2},
\; D(A_d)=\{x\in H^2(0,1):x(0)=x(1)=0\}.$ 
Then $A_d$ is self-adjoint and negative on $Z$. Let $\cX=Z\times Z$ and consider
$$Q=\begin{bmatrix} I&0\\ 0&I \end{bmatrix},\quad E=\begin{bmatrix} I&0\\ 0&0 \end{bmatrix}, \qquad A= \begin{bmatrix} A_d&I\\
-A_d&I \end{bmatrix}, \qquad D(A)=D(A_d)\times Z.$$
We have $\cX_1=\operatorname{ran}E=Z\times\{0\}$ and $\cX_2=\ker E=\{0\}\times Z$.
Let $D=\{x\in D(A):Ax\in\cX_1\}$. If $x=(x_1,x_2)^\top\in D$, then $x_2=A_dx_1$ and
$ Ax= \begin{bsmallmatrix} 2A_dx_1\\
0 \end{bsmallmatrix}$.
Hence 
\begin{equation*}
    \operatorname{Re}\langle Ax,x\rangle_{\cX}
=2\operatorname{Re}\langle A_dx_1,x_1\rangle_Z =-2\left\|\frac{dx_1}{d\xi}\right\|_Z^2\leq0,
\end{equation*}
so $A|_D$ is dissipative. Moreover, $x\in\cX_2\cap D$ implies $x_1=0$ and hence $x_2=A_dx_1=0$. Thus, $\cX_2\cap D=\{0\}$.

Let $(z_1,0)^\top\in\operatorname{ran}E$. Then 
$(E-A)x= \begin{bmatrix} z_1\\ 0 \end{bmatrix}$ reduces to $(I-2A_d)x_1=z_1$ and $x_2=A_dx_1$. Since $A_d$ is self-adjoint and negative, $I-2A_d:D(A_d)\to Z$ is boundedly invertible. Therefore,
$\operatorname{ran}E\subseteq\operatorname{ran}(E-A).$
The assumptions of Theorem \ref{thm1b} are thus satisfied, and the reduced operator is
$A_r=2A_d,\; D(A_r)=D(A_d),$
 generates a contraction semigroup on $\cX_1$.

We now show that the corresponding right $E$-resolvent estimate fails with respect to the original norm on $\cX$.
 For  $ \lambda >0 , $
$$(\lambda E-A)^{-1} E = \begin{bmatrix}
    (\lambda I-2 A_d )^{-1} & 0 \\ A_d (\lambda I-2 A_d )^{-1} & 0
\end{bmatrix}. $$
Let $\{e_k\}_{k\in\mathbb N}$ be normalized eigenfunctions of $A_d$, with $A_de_k=-\mu_ke_k$ and $\mu_k\to\infty$. For $x=(e_k,0)^\top$, we obtain
\begin{equation*}
(\lambda E-A)^{-1}Ex = \begin{bmatrix} \dfrac{1}{\lambda+2\mu_k}e_k\\ -\dfrac{\mu_k}{\lambda+2\mu_k}e_k
\end{bmatrix}.
\end{equation*}
Consequently, for every \(\lambda>0\),
\[
\begin{aligned}
\|(\lambda E-A)^{-1}E\|_{\mathcal L(X)}
&\geq
\left\|(\lambda E-A)^{-1}E
\begin{pmatrix}e_k\\0\end{pmatrix}\right\|_X\geq \frac{\mu_k}{\lambda+2\mu_k}.
\end{aligned}
\]
Letting \(k\to\infty\) yields
\[
\|(\lambda E-A)^{-1}E\|_{\mathcal L(X)}
\geq \frac12,\qquad \lambda>0.
\]

Thus, no constant $K>0$ can satisfy $\|(\lambda E -A)^{-1}E\|_{\mathcal L(\cX)} \leq\frac{K}{\lambda},
\; \lambda>0.$  Therefore, the right $E$-resolvent estimate fails with respect to the original state-space norm on $\cX$. 

\end{example}

 Example \ref{example-no-rad} shows that the assumptions ensuring contraction semigroup generation on the dynamic subspace are not, in general, sufficient to obtain $E$-radiality. For a specific class of DAEs, however, $E$-radiality can be recovered by imposing an additional condition on the algebraic component and by choosing the state spaces appropriately. This is illustrated in the following theorem.
 
\begin{theorem} \label{main-Th}
Assume the hypotheses of Theorem \ref{thm1b} with
$\cZ=\cX=\cX_1\oplus\cX_2$, $Q=I$, and
$E=\begin{bmatrix} I&0\\0&0\end{bmatrix}$. Define
$D:=\{x\in D(A):Ax\in\cX_1\}$. For $x\in\cX$, write
$x=x_1+x_2$, with $x_i\in\cX_i$. Suppose, in addition, that there exists
$A_a\in\mathcal L(\cX_1,\cX_2)$ such that $x_2=A_ax_1$ for every
$x=x_1+x_2\in D$. Define
$$\cX_c:=\{x_1+A_ax_1:x_1\in\cX_1\}, \; \text{and} \; \; \cZ_c:=E\cX_c=\cX_1\times\{0\}.$$
Then $\cX_c$ is a Hilbert space,
$A_D:=A|_D:D\subset\cX_c\to\cZ_c$ is closed and densely defined,
$E|_{\cX_c}\in\mathcal L(\cX_c,\cZ_c)$, and $A_D$ is $E$-radial with respect to the spaces 
$\cX_c$ and $\cZ_c$.
\end{theorem}

\begin{proof}
Since $A_a\in\mathcal L(\cX_1,\cX_2)$,
$\cX_c=\operatorname{graph}(A_a)$ is a closed subspace of
$\cX_1\oplus\cX_2$, and hence a Hilbert space. Moreover, for
$x=x_1+A_ax_1\in\cX_c$, we have $Ex=Ex_1\in\cZ_c$, so
$E|_{\cX_c}\in\mathcal L(\cX_c,\cZ_c)$.

By Theorem~\ref{thm1b}, the reduced operator
$A_r:D(A_r)\subset\cX_1\to\cX_1$ generates a contraction semigroup. Define $J:\cX_1\to\cX_c$ by $Jx_1=x_1+A_ax_1$. Since
$A_a\in\mathcal L(\cX_1,\cX_2)$, $J$ is bounded, and it is
surjective by the definition $\cX_c=J\cX_1$. Moreover,
$\cX_1\cap\cX_2=\{0\}$ implies that $J$ is injective. Since
$\cX_1\perp\cX_2$,
\[
\|Jx_1\|_{\cX}^2
=
\|x_1\|_{\cX}^2+\|A_ax_1\|_{\cX}^2
\geq \|x_1\|_{\cX}^2,
\]
and hence $\|J^{-1}y\|_{\cX}\leq \|y\|_{\cX}$ for every
$y\in\cX_c$. Therefore, $J$ is a bounded isomorphism. By the reconstruction assumption,
$D=J(D(A_r))$. Moreover,
\[
A_DJx_1=EA_rx_1,\qquad x_1\in D(A_r).
\]
Since $A_r$ is closed and densely defined and $J$ is an isomorphism,
$A_D$ is closed and densely defined.

For $s>0$ and $x_1\in D(A_r)$,
$(sE-A_D)Jx_1=E(sI-A_r)x_1$. Since $A_r$ generates a contraction
semigroup, $sI-A_r$ is boundedly invertible and
$\|(sI-A_r)^{-n}\|_{\mathcal L(\cX_1)}\leq s^{-n}$ for all
$n\in\mathbb N$. Hence $sE-A_D:D\to\cZ_c$ is boundedly invertible.

Since every $x\in\cX_c$ can be written uniquely as $x=Jx_1$ with
$x_1=J^{-1}x$, we obtain
$\bigl((sE-A_D)^{-1}E\bigr)^n
=J(sI-A_r)^{-n}J^{-1}$ on $\cX_c$, and therefore
\[
\left\|\bigl((sE-A_D)^{-1}E\bigr)^n\right\|_{\mathcal L(\cX_c)}
\leq \frac{\|J\|\,\|J^{-1}\|}{s^n}.
\]
Similarly,
$\bigl(E(sE-A_D)^{-1}\bigr)^n=(sI-A_r)^{-n}$ on $\cZ_c$, and hence
\[
\left\|\bigl(E(sE-A_D)^{-1}\bigr)^n\right\|_{\mathcal L(\cZ_c)}
\leq \frac{1}{s^n}.
\]
Thus both $E$-resolvent estimates hold, and $A_D$ is $E$-radial on
$\cX_c$ and $\cZ_c$.
\end{proof}

\subsection{Coupled Systems}

Let $Z$ be a complex Hilbert space. We consider the coupled system
\begin{subequations}\label{coupled_system}
\begin{align}
    \frac{d}{dt}p(t)&=A_1(p(t),y(t)),\\
    0&=A_2(p(t),y(t)), \label{alg-eq}
\end{align}
\end{subequations}
which includes models arising, for example, in poroelasticity
\cite{bociu2025analysis} and lithium-ion cells
\cite{chaturvedi2010algorithms}. {Here, $A_i:D(A_i)\subset Z\times Z\to Z$,
$i=1,2$, are closed linear operators.

Setting
\[
x=\begin{bmatrix}p\\y\end{bmatrix},\quad
E=\begin{bmatrix}I&0\\0&0\end{bmatrix},\quad A(p,y)=
\begin{bmatrix}
A_1(p,y)\\
A_2(p,y)
\end{bmatrix}, \] \color{black}
with $D(A)=
D(A_1)\cap D(A_2)$ assumed dense in $Z\times Z$,}
system \eqref{coupled_system} can be
written as
\begin{equation}
    \frac{d}{dt}Ex(t)=Ax(t).
\end{equation}

{The next result yields wellposedness and
$E$-radiality of the coupled system \eqref{coupled_system} } without requiring a
Schur-complement formulation, as used in
\cite[Section IV]{jacob2022solvability}.
{
\begin{proposition}\label{prop:coupled}
Define $
D:=\left\{
\begin{bmatrix}p\\y\end{bmatrix}\in D(A):
A_2(p,y)=0
\right\}$
and assume
\begin{enumerate}
    \item $A|_D$ is dissipative, that is,
    $
    \operatorname{Re}
    \left\langle
    A\begin{bmatrix}p\\y\end{bmatrix},
    \begin{bmatrix}p\\y\end{bmatrix}
    \right\rangle_{Z\times Z}
    \leq 0$ for $
    \begin{bmatrix}p\\y\end{bmatrix}\in D$.

    \item $Z\times\{0\}\subseteq\operatorname{ran}(E-A)$.

    \item There exists $A_a\in\mathcal{L}(Z)$ such that $y=A_ap$ for 
        $\begin{bmatrix}p\\y\end{bmatrix}\in D$.
\end{enumerate}
Then the reduced operator $A_r:D(A_r)\subset Z\to Z$,
defined by
\[
A_rp=A_1(p,A_a p),
\]
with domain
\[
D(A_r)=
\left\{
p\in Z:
\begin{bmatrix}p\\A_ap\end{bmatrix}\in D
\right\},
\]
generates a contraction semigroup on $Z$.
Moreover, $A_D:=A|_D$ is $E$-radial with respect to
the spaces
\[
X_c:=
\left\{
\begin{bmatrix}p\\A_ap\end{bmatrix}:p\in Z
\right\},
\qquad
Z_c:=Z\times\{0\}.
\]
\end{proposition}
}
\begin{proof}
We verify the hypotheses of Theorem~\ref{main-Th}.
Since $E=\left[\begin{smallmatrix}I&0\\0&0\end{smallmatrix}\right]$
is self-adjoint and nonnegative, and
$\operatorname{ran}E=Z\times\{0\}$ is closed,
the structural assumptions of Theorem~\ref{thm1b}
are satisfied with $Q=I$.
Moreover, $\cX_1=Z\times\{0\}$ and
$\cX_2=\{0\}\times Z$. 
The first two assumptions give, respectively, the
dissipativity of $A|_D$ and the range condition
$\operatorname{ran}E\subseteq\operatorname{ran}(E-A)$.

If $\left[\begin{smallmatrix}0\\y\end{smallmatrix}\right]\in \cX_2\cap D$,
then the third assumption implies $y=A_a0=0$.
Thus, $\cX_2\cap D=\{0\}$.

Finally, the third assumption ensures that the
algebraic component is determined by the dynamic
component through the bounded operator $A_a$.
Therefore, the hypotheses of Theorem~\ref{main-Th} are satisfied.
The conclusion follows from Theorems~\ref{thm1b} and~\ref{main-Th}.
\end{proof}

\begin{rem} {
Suppose that the  operator $A_2$ admits
the representation $A_2(p,y)=A_3p+A_4y$
for $(p,y)\in D(A)$.}
A sufficient condition for the third assumption
of the proposition is that $A_4$ is closed and bijective and
$A_4^{-1}A_3$ admits a bounded extension to $Z$.
Indeed, the algebraic constraint then yields
\[
y=-A_4^{-1}A_3p,
\]
and the bounded extension of $-A_4^{-1}A_3$
can be chosen as $A_a$.

Since $A_4$ is closed, its bijectivity already
implies $A_4^{-1}\in\mathcal{L}(Z)$ by the closed
graph theorem. However, the boundedness of
$A_4^{-1}A_3$ does not follow from the boundedness
of $A_4^{-1}$ alone.

However, even in this situation, where the operator $A_2$ decouples into $A_3p+A_4y , $ invertibility of $A_4$ is not necessary.
For example, let $Z=\ell^2(\mathbb N)$ and let
$S\in\mathcal L(Z)$ be the unilateral right shift.
Consider
\[
E=\begin{bmatrix}I&0\\0&0\end{bmatrix},
\qquad
A=\begin{bmatrix}-2I&I\\-S&S\end{bmatrix},
\qquad D(A)=Z\times Z.
\]
Here, $A_2(p,y)=S(y-p)$, so that
$A_3=-S$ and $A_4=S$.
Since $S$ is injective but not surjective,
$A_4=S$ is not bijective. However, the algebraic
constraint $S(y-p)=0$ implies $y=p$, so that
$A_a=I\in\mathcal L(Z)$. Moreover,
\[
D=\left\{\begin{bmatrix}p\\p\end{bmatrix}:p\in Z\right\},
\qquad
\operatorname{Re}\left\langle
A\begin{bmatrix}p\\p\end{bmatrix},
\begin{bmatrix}p\\p\end{bmatrix}
\right\rangle=-\|p\|^2\leq 0,
\]
and
\[
(E-A)\begin{bmatrix}p/2\\p/2\end{bmatrix}
=\begin{bmatrix}p\\0\end{bmatrix}.
\]
Thus, all assumptions of Proposition~\ref{prop:coupled}
are satisfied, and $A_r=-I$.
\end{rem}

\section{Semigroup Generation on Constrained Dynamic Subspaces} \label{section4}

In some DAEs, the algebraic equations impose additional constraints on the
dynamic variables, so that the admissible dynamics evolve on a closed
subspace of $\cX_1$. A related setting is considered in
\cite[Theorem~33]{zwart2024abstract}, where the dynamics are restricted
to a closed subspace of $\cX_1$ under a particular block structure for $E$.

Here, considering system \eqref{system1}, we extend this constrained-space framework without imposing a block
structure on $E,$  $Q$ and $A$. Under the structural assumptions specified below, we use the orthogonal decomposition $\cX=\cX_1\oplus\cX_2$ from
\eqref{orthogonal-decomposition}, where
$\cX_1=\ran(Q^\ast E)$ and $\cX_2=\ker E$, together with the notation
$\widetilde E:=Q^\ast E|_{\cX_1}$ and
$$D:=\{x\in D(AQ):Q^\ast AQx\in\cX_1\}.$$

\begin{theorem}\label{thm-constrained}
{Consider system \eqref{system1} and assume that
$\operatorname{ran}E$ is closed in $\cZ$ and
$E^\ast Q=Q^\ast E\geq 0$.}
Define
\[
W_0:=\{x_1\in\cX_1:\text{ there exists }x_2\in\cX_2
\text{ such that }x_1+x_2\in D\},
\qquad
\cX_0:=\overline{W_0}^{\,\cX_1}.
\]
Assume that
\begin{enumerate}
    \item $Q^\ast AQ|_D$ is linear and dissipative;
    \item $    \left\{  y_1\in\cX_0:     \text{ there exists }x_2\in\cX_2\cap D     \text{ such that }     Q^\ast AQx_2=\widetilde E y_1    \right\}    =\{0\};$ 
    \item $E\cX_0\subseteq\ran(E-AQ)$.
\end{enumerate}
Define $A_0:D(A_0)\subset\cX_0\to\cX_0$ by
$ A_0x_1:=\widetilde E^{-1}Q^\ast AQ(x_1+x_2),$ where
\[
D(A_0):= \left\{ x_1\in\cX_0: \begin{array}{l}
\text{there exists }x_2\in\cX_2\text{ such that }x_1+x_2\in D,\\
\widetilde E^{-1}Q^\ast AQ(x_1+x_2)\in\cX_0
\end{array}
\right\}.
\]
Then $A_0$ is well-defined and generates a contraction semigroup on
$\cX_0$, equipped with the inner product
$\langle x,y\rangle_{\cX_0}:=\langle x,\widetilde E y\rangle_{\cX}$.
\end{theorem}
\begin{proof}

Since $\cX_0$ is closed in $\cX_1$, it is a Hilbert space with the
inner product inherited from $\cX_1$.

We first show that $A_0$ is well-defined. Let $x_1\in D(A_0)$, and
suppose that $x_2,\widehat{x}_2\in\cX_2$ both satisfy the conditions
in the definition of $D(A_0)$. Since $D$ is a linear subspace,
$z_2:=x_2-\widehat{x}_2\in\cX_2\cap D$. Moreover, $y_1:=\widetilde E^{-1}Q^\ast AQz_2\in\cX_0.$  Thus $Q^\ast AQz_2=\widetilde E y_1$, and the second assumption implies that
$y_1=0$. Hence
$Q^\ast AQ(x_1+x_2)=Q^\ast AQ(x_1+\widehat{x}_2)$, so $A_0x_1$ is
independent of the choice of $x_2$. Since $D(A_0)$ is a linear subspace and the operators in the definition
of $A_0$ are linear, $A_0$ is linear.

Next, let $x_1\in D(A_0)$ and choose $x_2\in\cX_2$ as in the
definition of $D(A_0)$. Since $x_1+x_2\in D$ and
$x_2\perp\cX_1$,
\[
\begin{aligned}
\Re\langle A_0x_1,x_1\rangle_{\cX_0}
&=\Re\langle Q^\ast AQ(x_1+x_2),x_1\rangle_{\cX}=\Re\langle Q^\ast AQ(x_1+x_2),x_1+x_2\rangle_{\cX}
\leq0,
\end{aligned}
\]
where the last inequality follows from the first assumption. Thus $A_0$ is dissipative.

We now show that $\ran(I-A_0)=\cX_0$. Let $y_1\in\cX_0$. By
the third assumption, there exists $x=x_1+x_2\in D(AQ)$, with
$x_i\in\cX_i$, such that $Ey_1=(E-AQ)x$. Since
$x_2\in\ker E$, applying $Q^\ast$ gives
\[
Q^\ast AQ(x_1+x_2)=\widetilde E(x_1-y_1)\in\cX_1.
\]
Hence $x_1+x_2\in D$, so $x_1\in W_0\subseteq\cX_0$. Since also
$x_1-y_1\in\cX_0$, we have $x_1\in D(A_0)$ and
$A_0x_1=x_1-y_1$. Therefore $(I-A_0)x_1=y_1$, proving that
$\ran(I-A_0)=\cX_0$.

It remains to show that $D(A_0)$ is dense in $\cX_0$. Let
$y\in D(A_0)^\perp$. By the surjectivity of $I-A_0$, there exists
$x\in D(A_0)$ such that $y=(I-A_0)x$. Since
$\langle y,x\rangle_{\cX_0}=0$ and $A_0$ is dissipative,
\[
0
=\|x\|_{\cX_0}^2-\Re\langle A_0x,x\rangle_{\cX_0}
\geq\|x\|_{\cX_0}^2.
\]
Thus $x=0$ and hence $y=0$. Therefore
$D(A_0)^\perp=\{0\}$, so $D(A_0)$ is dense in $\cX_0$.

The operator $A_0$ is therefore densely defined and dissipative, with
$\ran(I-A_0)=\cX_0$. The Lumer-Phillips theorem now implies that
$A_0$ generates a contraction semigroup on $\cX_0$.
\end{proof}

Theorem \ref{thm-constrained} establishes contraction semigroup generation on the constrained dynamic subspace $\cX_0$. We now show that, when the algebraic component can be reconstructed from the dynamic component through a bounded operator, $E$-radiality can be obtained on suitable  spaces.

\begin{theorem}\label{thm-constrained-radial}

Assume the hypotheses of Theorem~\ref{thm-constrained}, and let
$A_0:D(A_0)\subset\cX_0\to\cX_0$ be the operator defined therein. Suppose, in addition, that there exists $A_a\in\mathcal L(\cX_0,\cX_2)$ such that $x_2=A_ax_1$ whenever $x_1\in\cX_0$, $x_2\in\cX_2$, and $x_1+x_2\in D$. Define $J:\cX_0\to\cX$ by $Jx_1:=x_1+A_ax_1$, and set \[
\cX_c:=J\cX_0 =\{x_1+A_ax_1:x_1\in\cX_0\},
\qquad
\cZ_c:=E\cX_c=E\cX_0.
\]
Define $A_c:D(A_c)\subset\cX_c\to\cZ_c$ by
\[
D(A_c):=J(D(A_0)),
\qquad
A_cJx_1:=AQJx_1,\quad x_1\in D(A_0).
\]
Then $A_c$ is $E$-radial with respect to the spaces $\cX_c$ and $\cZ_c$.
\end{theorem}

\begin{proof}
By Theorem~\ref{thm-constrained}, the operator
$A_0:D(A_0)\subset\cX_0\to\cX_0$ generates a contraction semigroup on
$\cX_0$. Hence, for every $s>0$, the operator
$sI-A_0:D(A_0)\to\cX_0$ is boundedly invertible and
$\|(sI-A_0)^{-n}\|_{\mathcal L(\cX_0)}\leq s^{-n}$ for all
$n\in\mathbb N$. As in the proof of Theorem \ref{main-Th}, the map $J:\cX_0\to\cX_c$, defined by
$Jx_1=x_1+A_ax_1$, is a bounded isomorphism. 
Since $A_a\cX_0\subseteq\cX_2=\ker E$, we have
$EJx_1=Ex_1$ for every $x_1\in\cX_0$, and hence
$\cZ_c=E\cX_c=E\cX_0$. The restriction
$E|_{\cX_0}:\cX_0\to\cZ_c$
is surjective by definition. Moreover, since
$\cX_0\subseteq\cX_1$ and $\widetilde E=Q^\ast E|_{\cX_1}$,
$\widetilde E x_1=Q^\ast Ex_1,
\; x_1\in\cX_0.$
Thus, for $z=Ex_1\in\cZ_c$,
$(E|_{\cX_0})^{-1}z
=\widetilde E^{-1}Q^\ast z.
$
Hence $E|_{\cX_0}:\cX_0\to\cZ_c$ is a bounded isomorphism.

We next relate $A_c$ to the reduced operator $A_0$. Let $x_1\in D(A_0)$. By the definition of $D(A_0)$, there exists $x_2\in\cX_2$ such that $x_1+x_2\in D$. By the reconstruction assumption, $x_2=A_ax_1$, and hence $Jx_1=x_1+A_ax_1\in D$. By the definition of $A_0$, we have
\[
Q^*AQJx_1=\widetilde E A_0x_1
          =Q^*EA_0x_1.
\]
Since $Q^*$ is injective, it follows that
\[
AQJx_1=EA_0x_1\in EX_0=Z_c.
\]
Therefore,
\begin{equation}\label{Ac-A0-relation}
A_cJx_1=EA_0x_1,\qquad x_1\in D(A_0).
\end{equation}
 Since $D(A_c)=J(D(A_0))$, the density of $D(A_0)$ in $\cX_0$ and the fact that $J$ is an isomorphism imply that $D(A_c)$ is dense in $\cX_c$. Moreover, \eqref{Ac-A0-relation} gives $A_c=(E|_{\cX_0})A_0J^{-1}$. Since $A_0$ is closed and $J$ and $E|_{\cX_0}$ are bounded isomorphisms, $A_c$ is closed.

It remains to verify the $E$-resolvent estimates. For $x_1\in D(A_0)$, using $EJx_1=Ex_1$ and \eqref{Ac-A0-relation}, we obtain $(sE-A_c)Jx_1=E(sI-A_0)x_1$. Therefore, for every $s>0$,
\begin{equation}
(sE-A_c)^{-1}
=
J(sI-A_0)^{-1}(E|_{\cX_0})^{-1}.
\label{resolvent-inverse}
\end{equation}

If $x=Jx_1\in\cX_c$, then $(E|_{\cX_0})^{-1}Ex=x_1=J^{-1}x$. Hence $(sE-A_c)^{-1}E=J(sI-A_0)^{-1}J^{-1}$ on $\cX_c$, and therefore $\bigl((sE-A_c)^{-1}E\bigr)^n=J(sI-A_0)^{-n}J^{-1}$. Thus,
\[
\left\|\bigl((sE-A_c)^{-1}E\bigr)^n\right\|_{\mathcal L(\cX_c)}
\leq \frac{\|J\|\,\|J^{-1}\|}{s^n}.
\]

Similarly, $E(sE-A_c)^{-1}=(E|_{\cX_0})(sI-A_0)^{-1}(E|_{\cX_0})^{-1}$ on $\cZ_c$, and hence
\[
\left\|\bigl(E(sE-A_c)^{-1}\bigr)^n\right\|_{\mathcal L(\cZ_c)}
\leq
\frac{\|E|_{\cX_0}\|\,\|(E|_{\cX_0})^{-1}\|}{s^n}.
\]

Setting $K:=\max\{\|J\|\,\|J^{-1}\|,\|E|_{\cX_0}\|\,\|(E|_{\cX_0})^{-1}\|\}$, both $E$-resolvent estimates hold for every $s>0$ and $n\in\mathbb N$. Therefore, $A_c$ is $E$-radial on $\cX_c$ and $\cZ_c$.
\end{proof}

The linearized Navier-Stokes system is shown in \cite[Example~36]{zwart2024abstract} to satisfy the assumption of Theorem \ref{thm-constrained}. In \cite[sect. 5.2]{sviridyuk2003linear} it is shown, by reformulation of the equations, that this system is radial of degree 1. This suggests that this class of systems, unlike those addressed by Theorem \ref{thm-constrained-radial}, are radial of higher degree, when radiality does hold.

\section{Conclusion} 

A wellposedness framework for a class of irregular infinite-dimensional PDAEs was described in which the operator pencil need not be surjective on the full state space. By combining dissipativity with a range condition restricted to the dynamic component, we obtained contraction semigroup generation on the corresponding dynamic subspace under assumptions weaker than those  imposed in \cite{jacob2022solvability,zwart2024abstract,sviridyuk2003linear}.
The presentation also provides a unified treatment of these results. 

 We further extended this framework to systems in which the algebraic equations impose additional constraints on the dynamic variables. In this case, the dynamics evolve on a proper closed subspace $\cX_0$, on which the reduced operator generates a contraction semigroup. 
 
We also examined the relation between this reduced semigroup formulation and $E$-radiality. The wellposedness assumptions alone are not sufficient to guarantee $E$-radiality on the original state space, as demonstrated by Example~\ref{example-no-rad}. However, when the algebraic component can be reconstructed from the dynamic component through a bounded operator, the system admits an $E$-radial formulation on suitably chosen state and range spaces.

Several questions remain open. In particular, it would be useful to identify conditions leading to higher $(E,p)$-radiality and also to determine whether the bounded reconstruction assumption can be weakened, for example by allowing the reconstruction operator to be unbounded or defined only on the domain of the reduced operator. 
 
\bibliographystyle{ieeetr}     
\bibliography{ref}  
\end{document}